\documentclass[12pt]{amsart}
\usepackage{amssymb}
\usepackage{amsmath}
\usepackage{pictex}
\newcommand\AND{\quad\text{and}\quad}
\newcommand\Aut{\operatorname{\sf Aut}}
\newcommand\BD{\operatorname{\sf B}}
\newcommand\bd{\partial}
\newcommand\CC{\mathcal C}

\newcommand\etab{\boldsymbol{\eta}}

\newcommand\Ga{\mathcal G}
\newcommand\geo[1]{\overline{#1}}
\newcommand\gb{\mathbf{g}} 
\newcommand\id{\text{\sl id}}
\newcommand\K{\mathbb K}
\newcommand\Kr{\K_r}
\newcommand\N{\mathbb N}
\newcommand\Prob{\mathsf{Pr}}
\newcommand\R{\mathbb R}
\newcommand\supp{\operatorname{\sf supp}}
\newcommand\T{\mathbb T}
\newcommand\uf{\mathfrak{u}}

\newcommand\vf{\mathfrak{v}}

\newcommand\wh{\widehat}
\newcommand\whCC{\,\widehat{\!\CC}}
\newcommand\wt{\widetilde}
\newcommand\xb{\mathbf{x}} 
\newcommand\Z{\mathbb Z}
\newcommand\zero{\text{\sl 0}}
\newcommand\Zr{{\mathcal Z}_r}

\numberwithin{equation}{section}

\newtheoremstyle{mythm}
  {9pt}
  {9pt}
  {\itshape}
  {0pt}
  {\bfseries}
  {}
  { }
  {\thmnumber{(#2)}\thmname{ #1}\thmnote{ #3}}

\newtheoremstyle{mydef}
  {9pt}
  {9pt}
  {\normalfont}
  {0pt}
  {\bfseries}
  {}
  { }
  {\thmnumber{(#2)}\thmname{ #1}\thmnote{ #3}}

\theoremstyle{mythm}
\newtheorem{thm}[equation]{Theorem.}
\newtheorem{pro}[equation]{Proposition.}
\newtheorem{lem}[equation]{Lemma.}

\theoremstyle{mydef}

\begin{document}$\,$ \vspace{-1truecm}
\title{The Poisson boundary of \\ lamplighter random walks on trees}
\author{\bf Anders KARLSSON and Wolfgang WOESS}
\address{\parbox{.8\linewidth}{Department of Mathematics,
Royal Institute of Technology\\
100 44 Stockholm, Sweden\\}}
\email{akarl@math.kth.se}
\address{\parbox{.8\linewidth}{Institut f\"ur Mathematik C, 
Technische Universit\"at Graz,\\
Steyrergasse 30, A-8010 Graz, Austria\\}}
\email{woess@TUGraz.at}
\date{June 1, 2006} 
\thanks{Supported by ESF program RDSES and 
by Austrian Science Fund (FWF) P15577}
\begin{abstract}
Let $\T_q$ be the homogeneous tree with degree $q+1 \ge 3$ and $\Ga$ a finitely
generated group whose Cayley graph is $\T_q$. The associated lamplighter group is
the wreath product $\Zr \wr \Ga$, where $\Zr$ is the cyclic group of order $r$. 
For a large class of random walks on this group,
we prove almost sure convergence to a natural geometric boundary. If the probability 
law governing the random walk has finite first moment, then the probability space
formed by this geometric boundary together with the limit distribution of the
random walk is proved to be maximal, that is, the Poisson boundary.
We also prove that the Dirichlet problem at infinity is solvable 
for continuous functions on the active part of the boundary, if the lamplighter
``operates at bounded range''.
 \end{abstract}

\maketitle

\markboth{{\sf A. Karlsson and W. Woess}}
{{\sf Lamplighter random walks on trees}}
\baselineskip 15pt

\section{Introduction}\label{intro}

Let $\T = \T_q$ be the homogeneous tree with degree $q+1\geq 3$. Assume that at each vertex $x \in \T$
there is a lamp which may be switched off or on in $r$ different states
of ``intensity'', encoded by the set
$\{0,\dots, r-1\}$, where the state $0$ represents ``off''. 
We think of $\{0,\dots, r-1\}$ as the vertex set of the complete graph
$\Kr$ (all pairs of distinct elements are adjacent). 
As an introductory example,
consider the following random process: at the beginning, all lamps are switched
of. A ``lamplighter'' starts a random walk at a vertex of $\T$. With 
probability $\theta$, he chooses to move, that is, he makes a step to a 
randomly selected neighbouring vertex (without changing the lamps). With 
probability $1-\theta$, he chooses to ``switch'', that is, he randomly modifies the state 
of the lamp where he stands (and does not move). 
At each step we observe the actual position in the tree and the 
configuration of the lamps that are switched on. Thus, our process evolves
on the state state space $\Kr \wr \T$ consisting of pairs $(\eta,x)$, 
where $x \in \T$ and the configuration $\eta: \T \to \Kr$ is a 
function with finite support. The set of all configurations 
is denoted $\CC$.  It can be equipped with the structure of an abelian
group with pointwise addition modulo $r$, and writing $\eta - \eta'$ below
refers to this operation.
In the sequel, we shall often write $\eta\,x$ instead of $(\eta,x)$.

We equip $\Kr \wr \T$ with a 
graph structure, where the neighbourhood relation is given by
\begin{equation}\label{eq:nbhd}
(\eta,x) \sim (\eta',x') : \iff \begin{cases}
x \sim x' \;\text{in}\; \T \;\text{and}\; \eta = \eta'\,, \quad\text{or}&\\
x = x' \;\text{in}\; \T \;\text{and}\; \supp(\eta - \eta') = \{x\}\,.&
\end{cases}
\end{equation}
We view $\T$ as a Cayley graph of a finitely generated group $\Ga$, 
a free product of two-element and infinite cyclic groups. Thus, vertices of $\T$
are (in one-to-one correspondence with) elements of $\Ga$, and $x^{-1}$ and
$xx'$ refer to the operations in this group. By a slight deviation from the usual
notation, we write $o$ for the group identity and think of it as a root (origin)
of the tree.  
Then our \emph{lamplighter graph} $\Kr \wr \T$ is a Cayley graph of the 
\emph{wreath product} $\Gamma = \Zr \wr \Ga$, where $\Zr = \Z/(r\Z)$. More precisely, 
every $x \in \Ga$ acts on $\CC$ by the translation $T_x$, where 
$T_x\eta(y) = \eta(x^{-1}y)$. The resulting semidirect product is
\begin{equation}\label{eq:semidir}
\Zr \wr \Ga = \CC \rtimes \Ga\,,\quad \text{with} \quad
(\eta,x)(\eta',x') = (\eta + T_{x}\eta', xx')\,,
\end{equation}
Our random process is a Markov chain $Z_n=(Y_n,X_n)$ on $\Kr \wr \T$, where
$X_n$ is the position and $Y_n$ the configuration at time $n$. Its 
one-step transition probabilities 
$p(\eta\,x,\eta'\,x') = \Prob[Z_{n+1} = \eta'\,x' | Z_n = \eta\,x]$ are given
by
\begin{equation}\label{eq:transprob}
p(\eta\,x,\eta'\,x') = \begin{cases} \theta/(q+1)\,,&\text{if}\; x \sim x'\;
                       \text{and}\; \eta = \eta'\,,\\
		       (1-\theta)/(r-1)\,,&\text{if}\; x = x'
                       \;\text{and}\;\supp(\eta - \eta') = \{ x \}\,\\
                                     0 &\text{in all other cases.}
		       \end{cases}
\end{equation}
Thus, $p(\eta\,x,\eta'\,x') = \mu\bigl( (\eta\,x)^{-1}(\eta'\,x')\bigr)$,  
where $\mu$ is a probability measure on the group $\Zr \wr \Ga$. 
In particular, we can view $Z_n$ as a random walk on that group.
The random walk is transient, that is, with probability $1$, it visits every
finite subset of $\Kr \wr \T$ only finitely many times. Thus, $Z_n$ tends
``to infinity'', and the purpose of this note is to relate this property 
in a more detailed way with the underlying structure. 
Below we shall also consider more general random walks on $\Zr \wr \Ga$. 

We shall determine the \emph{Poisson boundary} of a general class of 
lampligher random walks on $\Kr \wr \T$, resp. $\Zr \wr \Ga$ that includes 
the basic example  \eqref{eq:transprob}. This boundary can be defined in several 
equivalent ways, see {\sc Kaimanovich and Vershik~\cite{KaVe}} and
{\sc Kaimanovich~\cite{Kai2}}. The Poisson boundary of a random walk
on a group is a measure space, determined uniquely up to
measure-theoretic isomorphism. One quick definition
is to say that it is the space of ergodic components in the trajectory space
of the random walk. Another approach is via bounded harmonic functions, see
below. Here, we take a more topological viewpoint. We attach to
$\Kr \wr \T$ a natural boundary $\Pi$ at infinity, defined purely in geometric
terms, such that $(\Kr \wr \T) \cup \Pi$ is a metrizable
space (not necessarily compact or complete) on which the group $\Gamma=\Zr \wr \Ga$ 
acts by homeomorphisms, and every point in $\Pi$ is the accumulation
point of a sequence in $\Kr \wr \T$. We then show that in that topology,
$(Z_n)$ converges almost surely to a $\Pi$-valued random variable 
$Z_{\infty}$. Let $\nu$ be the distribution of $Z_{\infty}$, given
that the initial position and configuration of the random walk at time
$n=0$ are $o$ and the zero configuration $\zero$. The measure $\nu$ is often called the
\emph{harmonic measure} or \emph{limit distribution}. The pair 
$(\Pi,\nu)$ provides a model for the behaviour at infinity 
(in time and space) of the random walk. We give a quite simple proof 
that this is indeed
the Poisson boundary under rather general assumptions. This tells us that
up to sets with measure $0$, we have found the finest model for
distinguishing limit points at infinity of the random walk.
The (geometric) tool that we shall use for proving that $(\Pi,\nu)$ 
is the Poisson boundary is the \emph{strip criterion} of 
{\sc Kaimanovich~\cite[\S 6]{Kai2}, \cite[Thm. 5.19]{KaWo}}.

In section \ref{Dirichlet}, we prove that the Dirichlet 
problem is solvable with respect
to this natural geometric boundary. Again, in the spirit of this 
article, the focus is not on proving the most general result
possible.

Let us conclude the introduction wih a brief and incomplete overview
of recent work on lamplighter random walks and identifications of the 
Poisson boundary.

The abovementioned paper of {\sc Kaimanovich and Vershik~\cite{KaVe}}
may serve as a major source for the previous literature, different
equivalent definitions of the Poisson boundary, and a wealth of results
and methods. There one also finds the first interesting results
on the Poisson boundary of lamplighter random walks, namely on $\Zr \wr \Z^d$.
The techniques were refined in the subsequent body of work of
{\sc Kaimanovich,} see e.g. \cite{Kai1}, \cite{Kai2}, \cite{KaWo} and the 
references therein. Within the study of random walks on groups, wreath products
(lamplighter walks) have been the object of intensive studies in the last
decade. Wreath products exhibit interesting types of asymptotic behaviour
of $n$-step return probabilities, see {\sc Saloff-Coste and
Pittet}~\cite{PitSal1}, \cite{PitSal2}, {\sc Revelle}~\cite{Rev2}. 
The rate of escape has been studied by 
{\sc Lyons, Pemantle and Peres}~\cite{LyoPemPer}, {\sc Erschler}~\cite{Er1}, 
\cite{Er2}, {\sc Revelle}~\cite{Rev1} and, for simple lamplighter walks on 
trees, by {\sc Gilch}~\cite{Gi}. For the spectrum of transition operators, see
{\sc Grigorchuk and \.Zuk}~\cite{GrZu},  {\sc Dicks and Schick}~\cite{DiSc}
and  {\sc Bartholdi and Woess}~\cite{BaWo}. The positive harmonic functions
for a class of random walks on \hbox{$\Zr \wr \Z$} have been determined
by {\sc Woess}~\cite{Wo1} and {\sc Brofferio and Woess}~\cite{BrWo2}, who have
also determined the full Martin compactification in a ``nearest neighbour'' 
case~\cite{BrWo1}. For the Poisson boundary and bounded harmonic functions
for various types of lamplighter random walks,
besides \cite{KaVe} and \cite{Kai1}, see once more the  
impressive work of {\sc Erschler}~\cite{Er3}.

Concerning the Dirichlet problem in the discrete setting we refer to 
{\sc Woess}~\cite[Chapter IV]{Wbook} for more information.

\section{Convergence to the geometric boundary}\label{convergence}

The lamplighter graph $\Kr \wr \T$ with neighbourhood defined \eqref{eq:nbhd}
is far from being tree-like (it has one end). 
It is easy (and well known) to describe the graph metric. A shortest path from
$\eta\,x$ to $\eta'\,x'$ in the lamplighter graph must be such that the
lamplighter starts at $x$, walks along the tree and visits every 
$y \in \supp(\eta' - \eta)$, where he has to switch the lamp from state
$\eta(y)$ to state $\eta(y')$, and at the end, he must reach $x'$. 
Thus, $d(\eta\,x,\eta'\,x') = \ell + |\supp(\eta' - \eta)|$, where
$\ell$ is the smallest length of a ``travelling salesman'' tour (walk)
from $x$ to $x'$ that visits each element of $\supp(\eta' - \eta)$.
This description of the metric does not require that the base graph is a 
tree, however note that in a tree, the ``travelling salesman'' algorithm required 
for finding such a tour is easy to implement, see e.g. {\sc Parry~\cite{Pa}}
and {\sc Ceccherini-Silberstein and Woess~\cite[Example 2]{CeWo}}.

Since the group $\Gamma=\Zr \wr \Ga$ is non-amenable, the random walks 
$Z_n = (Y_n,X_n)$ that we consider here are all \emph{transient}, that is, 
$d(Z_n,Z_0) \to \infty$  almost surely, see \cite[Thm. 3.24]{Wbook} for this
result going back to {\sc Kesten~\cite{Ke1}, \cite{Ke2}}. 
The main question considered here is whether we can describe 
in a more detailed, geometric way how $(Z_n)$ behaves at infinity. 

For this purpose, we first briefly recall the \emph{end compactification} of $\T$, 
whose graph metric we denote also by $d(\cdot,\cdot)$. 
A \emph{geodesic path}, resp. \emph{geodesic ray}, resp. \emph{infinite
geodesic} in $\T$ is a finite, 
resp. one-sided infinite, resp. doubly infinite sequence $(x_n)$ of vertices
of $\T$ such that $d(x_i,x_j) = |i-j|$ for all $i, j$. 
Two rays are \emph{equivalent} if, as sets, their symmetric difference is finite.
An \emph{end} of $\T$ is an equivalence class of rays. The space of 
ends is denoted $\bd \T$, and we write $\wh \T = \T \cup \bd \T$. 
For all $w, z \in \wh \T$ there is a unique geodesic $\geo{w\,z}$ 
that connects the two. In particular, if $x \in \T$ and $\uf \in \bd \T$ then 
$\geo{x\,\uf}$ is the ray that starts at $x$ and represents $\uf$.
Furthermore, if $\uf, \vf \in \bd \T$ ($\uf \ne \vf$) then
$\geo{\uf\,\vf}$ is the infinite geodesic whose two halves (split at any
vertex) are rays that respresent $\uf$ and $\vf$, respectively.
If $w, z \in \wh \T$, then their \emph{confluent} 
$c=w \wedge z$ with respect to the root vertex $o$ (=the identity element of $\Ga$) 
is defined by $\geo{o\,w} \cap \geo{o\,z} = \geo{o\,c}$. Let $(w|z) = d(o,c)$, which 
is finite unless $w=z \in \bd\T$. We can define a metric $\rho$ on $\wh \T$ by
\begin{equation}\label{eq:metric}
\rho(w,z) = \begin{cases} q^{-(w|z)}\,,&\text{if}\; z \ne w\,,\\
                          0\,,&\text{if}\;z=w\,.
            \end{cases}
\end{equation}
This makes $\wh \T$ a a compact ultrametric space with $\T$ as a dense, discrete 
subset. Each isometry $g \in \Aut(\T)$ extends to a homeomorphism of $\wh\T$.

The natural compactification of $\CC$ in the topology of pointwise convergence
is the set $\,\wh{\!\CC} = \Kr^{\T} = \{ \zeta : \T \to \K_r \}$ of \emph{all},
finitely or infinitely supported, configurations. Since the vertex set of the
lamplighter graph is $\CC \times \T$, the space 
$\wh{\Kr \wr \T} = \whCC \times \wh\T$ is a natural compactification, and
$\bd(\Kr \wr \T) = (\whCC \times \wh\T) \setminus (\CC \times \T)$ is a natural
``geometric'' boundary at infinity of the lamplighter graph. We shall see that this
boundary contains many points where the random walk $(Z_n)$ does not accumulate.
We define 
\begin{equation}\label{eq:Pi}
\begin{aligned}
\Pi &= \bigcup_{\uf \in \bd \T} \CC_{\uf} \times \{\uf\}\,,\quad\text{where}\\ 
\CC_{\uf} &= \{\zeta \in \whCC : 
\supp(\zeta) \;\text{is finite or accumulates just at}\; \uf\}\,.
\end{aligned}
\end{equation}
Since $\CC_{\uf}$ is dense in $\whCC$, the closure of $\Pi$ is the part 
$\whCC \times \bd\T$ of $\bd(\Kr \wr \T)$.
The action of the group $\Zr \wr \Ga$ on $\Kr \wr \T$ by multiplication from the
left extends by homeomorphisms to $\wh{\Kr \wr \T}$ and leaves the Borel subset 
$\Pi$ invariant. Indeed, if $g=(\eta,x) \in \Gamma=\Zr \wr \Ga$ and 
$\beta = (\zeta,\uf) \in \whCC \times \wh\T$ then, precisely as in \eqref{eq:semidir},
\begin{equation}\label{eq:action}
g\beta = (\eta,x)(\zeta,\uf) = (\eta + T_x\zeta, x\uf)\,.
\end{equation}
(Addition of configurations is pointwise modulo $r$.)
If in addition $\zeta \in \CC_{\uf}$, where $\uf \in \bd \T$, then 
$\eta +  T_x\zeta \in \CC_{x\uf}$, since adding $\eta$ modifies
$T_x\zeta$ only in finitely many points.  

For the basic example \eqref{eq:transprob}, it is quite clear that $Z_n$ converges
almost surely to a random element $Z_{\infty} = (Y_{\infty},X_{\infty}) \in \Pi$.
Indeed, the $\T$-coordinate $X_n$ of $Z_n$ is the random walk on $\T$
with transition probabilities $\wt p(x,x') = \theta/(q+1)$, if $x' \sim x$, and
$\wt p(x,x) = 1-\theta$ (and $\wt p(x,y)=0$ if $d(y,x) \ge 2$). An elementary argument
using only transience yields that $X_n$ converges to a random elment 
$X_{\infty}\in \bd\T$. Also, only the states of those lamps can be modified which are
visited by $(X_n)$, and by transience, each vertex is visited finitely often: after
the last visit at a given vertex, the state of lamp sitting there remains unchanged.
Therefore $Y_n$ must converge pointwise to a random configuration $Y_{\infty}$
which can accumulate at no point besides $X_{\infty}$. 

We shall prove an analogous result for a much larger class of random walks on
$\Gamma=\Zr \wr \Ga$, or equivalently, on $\Kr \wr \T$. They are specified by a 
probability measure $\mu$ on the group $\Gamma$, and we suppose that 
$\supp(\mu)$ generates $\Gamma$.  
We can model the random walk on the probability space $(\Omega,\Prob)$,
where $\Omega = \Gamma^{\N}$ and $\Prob=\mu^{\N}$, with $\N=\{1,2,3,\ldots\}$.
The $n$-th projections $\gb_n=(\etab_n, \xb_n): \Omega \to \Z_r \wr \Ga$ are a sequence of 
independent, $\Gamma$-valued random variables with common distribution $\mu$. If 
$g_0 = (\eta_0,x_0) \in \Gamma$ then the sequence of random variables
$$
Z_n = (Y_n,X_n) = g_0 \gb_1 \cdots \gb_n\,,\; n \ge 0
$$
is the \emph{right random walk} on $\Gamma$ with law $\mu$ and starting point $g_0$. 
Its one-step transition probabilities are given by $p(g,g') = \mu( g^{-1}g')$, 
where $g= \eta\,x$ and $g' = \eta'\,x' \in \Gamma$.  Note that
\begin{equation}\label{eq:xnyn}
X_n = x_0 \xb_1 \cdots \xb_n \AND Y_n = \eta_0 + \sum_{k=1}^n T_{X_{k-1}}\etab_k\,,
\end{equation}
with pointwise addition modulo $r$. In particular, $(X_n)$
is the right random walk on the ``base'' group $\Ga$ with starting point $x_0$
whose law ist the projection $\wt \mu$ of $\mu$,
\begin{equation}\label{eq:mutilde}
\wt \mu (x) = \sum_{\eta} \mu(\eta\,x)\,.
\end{equation}
Before stating a result on convergence of $(Z_n)$ to $\Pi$, we have to specify additional
properties. We say that the lamplighter \emph{operates at bounded range,} if
\begin{equation}\label{eq:lamprange}
R(\mu) = \max\bigl\{ \min\{d(y,o), d(y,x)\} : 
\mu(\eta\,x) > 0\,,\; y  \in \supp(\eta) \bigr\} < \infty\,.
\end{equation}
This means that when the lamplighter steps from $x$ to $x'$ in $\T$ then while 
doing so, he can modify the current configuration of lamps only at points which are at 
bounded distance from $x$ or $x'$. (Note that here wo do not require that
$\mu$ itself has finite support.)
  
The law $\mu$ of the random walk is said to have \emph{finite first moment,} if 
\begin{equation}\label{eq:moment}
\sum_{\eta\,x \in \Ga} d(\zero\,o,\eta\,x)\,\mu(\eta\,x) < \infty\,.
\end{equation}
In this case, it is a well known consequence of Kingman's subbadditive ergodic 
theorem (see e.g. {\sc Derriennic~\cite{De}}) that there are finite constants 
$\ell(P)$ and $m(P)$ such that
$$
\frac{d(Z_n,Z_0)}{n} \to \ell(P) \AND  \frac{d(X_n,X_0)}{n} \to m(P) \quad
\text{almost surely.}
$$
It is clear that $\ell(P) \ge m(P)$. Recent \cite{Gi} and ongoing work by {\sc Gilch}
suggests that $\ell(P) > m(P)$ strictly.  
\begin{lem}\label{lem:rate}
If $P$ has finite first moment then $m(P) > 0$.
\end{lem}

\begin{proof} This follows from the fact that $\Ga$ is non-amenable, see 
\cite[Thm. 8.14]{Wbook}. Again, for random walks with finite range on a group,
this goes back to {\sc Kesten~\cite{Ke1}, \cite{Ke2}}.
\end{proof}

\begin{thm}\label{thm:conv} Let $Z_n=(Y_n,X_n)$ be a random walk with law
$\mu$ on $\Gamma =\Zr \wr \Ga$, such that $\supp(\mu)$ generates $\Gamma$.

If the lamplighter operates at finite range, or if $\mu$ has finite first moment,
then there is a $\Pi$-valued random variable $Z_{\infty} = (Y_{\infty},X_{\infty})$ 
such that $Z_n \to Z_{\infty}$ in the topology of $\wh{\Kr \wr \T}$, almost surely
for every starting point $g_0=\eta_0\,x_0$. 

Furthermore, the distribution of $X_{\infty}$ is a continuous measure on $\bd \T$
(it carries no point mass), and consequently the same is true for the distribution
of $Z_{\infty}$ on $\Pi$.
\end{thm}

\begin{proof} We may suppose without loss of generality that $g_0=\id$,
where $\id = (\zero,o)$ is the identity of $\Zr \wr \Ga$.

The law $\wt\mu$ of the projected random walk $(X_n)$ on $\Ga$ is such that
its support generates $\Ga$. 
{\sc Cartwright and Soardi~\cite{CaSo}} have shown that without any moment assumption,
such a random walk on $\T\equiv\Ga$ converges almost surely to a random end,
that is, a $\bd\T$-valued random variable $X_{\infty}$. 

Now suppose first that the lamplighter operates at bounded range. 
Let $(y_n)$ be an unbounded sequence in $\T$ with $y_n \in \supp(Y_n)$, i.e., 
$y_n$ is a vertex where the lamp is ``on'' at time $n$. Then 
we see from \eqref{eq:xnyn} that $y_n$ must be at bounded distance from the
initial trajectory $\{ X_0, X_1, \dots, X_n\}$.
Therefore, we must have $y_n \to X_{\infty}$. Indeed, one sees immediately from
the definition \eqref{eq:metric} of the metric $\rho$ that the end 
compactification $\wh\T$  has the following property: if $(x_n)$, $(y_n)$ are
two sequences in $\T$ such that $x_n \to \uf  \in \bd \T$ and $\sup_n d(y_n,x_n)<\infty$,
then also $y_n \to \uf$.

Next, assume that the random walk has finite first moment.
Define the integer-valued random variables $M_n = \max \{ d(y,o) : y \in \supp(\etab_n) \}$.
They are independent and identically distributed and also have finite first moment.
Therefore $M_n/n \to 0$ almost surely. This implies that the following holds with 
probability $1$.
\begin{equation}\label{eq:conv}
\begin{gathered}
\text{If $(y_n)$ is a sequence in $\T$ such that 
$y_n \in \supp\bigl(T_{X_{n-1}}\etab_n\bigr)$ 
for each $n$, then}\\
d(y_n,X_{n-1})/n \to 0 \quad \text{as}\quad n \to \infty\,.
\end{gathered}
\end{equation}
Recall that almost surely $X_n  \to X_{\infty} \in \bd\T$ and $d(X_n,o)/n \to m(P)$, with
$m(P) > 0$ by Lemma \ref{lem:rate}. Thus, we get for the confluents 
$c_n = y_n \wedge X_{n-1}$ that $d(c_n,X_{n-1})/n \le d(y_n,X_{n-1})/n \to 0$, whence
$d(c_n,o)/n \to m(P)$. Therefore, for the metric \eqref{eq:metric} of $\wh\T$, we
obtain $\rho(y_n,X_{n-1}) \to 0$. Consequently $y_n \to X_{\infty}$.

Now observe that by \eqref{eq:xnyn},
$$
\supp(Y_n) \subset \bigcup_{k=1}^n \supp\bigl(T_{X_{k-1}}\etab_k\bigr)\,,
$$
which is a union of finite sets. Thus, almost surely by the above, 
if $(y_n)$ is an unbounded sequence in $\T$ with $y_n \in \supp(Y_n)$, 
it must converge to $X_{\infty}$.  

Finally, from \cite{CaSo} it is also known (without moment hypothesis)
that the distribution of $X_{\infty}$ is
a continuous measure on $\bd\T$. 
\end{proof}

\section{The Poisson boundary}\label{Poisson}

Under the assumptions of Theorem \ref{thm:conv}, we can define the limit distribution
$\nu$ of the random walk. This is the probability measure defined for Borel sets
$U \subset \Pi$ by
\begin{equation}\label{eq:nu} 
\nu(U) = \Prob[Z_{\infty}\in U \mid Z_0 = \id]
\end{equation}
Recall the natural action \eqref{eq:action} of $\Gamma$ on $\Pi$.
If $g_0 \in \Gamma$, then 
$$
\Prob[Z_{\infty}\in U \mid Z_0 = g] = \nu(g^{-1}U)\,.
$$
This implies that $\nu$ satisfies the convolution equation $\mu*\nu = \nu$.
In particular, the measure space $(\Pi,\nu)$ is a \emph{boundary} of the
random walk on $\Gamma$ in the sense of {\sc Furstenberg~\cite{Fu}}.

A general boundary is a suitable probability space $(\BD,\lambda)$ such that
$\Gamma$ acts on $\BD$ by measurable bijections and $\mu*\lambda=\lambda$.
As is explained in \cite{Fu}, one can then typically (i.e., when $\BD$
carries a topology and $\Gamma$ acts continuously) construct a topology on 
$\Gamma \cup \BD$
such that $(Z_n)$ converges almost surely to a $\BD$-valued random variable
whose distribution is $\lambda$, given that $Z_0=\id$. This means that 
$(\BD,\lambda)$ is a model for describing in detail how $(Z_n)$ tends to
infinity, that is, for distinguishing limit points of $(Z_n)$ as $n \to \infty$.
When comparing two boundaries, this has of course to be done modulo sets
with measure $0$. For example, in the case of the boundary $(\Pi,\nu)$ it may appear more
natural to consider $\nu$ as a measure on the closure 
$\overline{\Pi} = \whCC \times \bd\T$ of
$\Pi$ in $\wh{\Kr \wr \T}$, with the same definition as in \eqref{eq:nu} but
Borel sets $U \subset \whCC \times \bd\T$. Since $\nu$ charges only the (dense)
subset $\Pi$, the measure spaces $(\Pi,\nu)$ and $(\whCC \times \bd\T, \nu)$ are
isomorphic.

In this spirit, our question is whether the boundary $(\Pi,\nu)$ is \emph{maximal.}
This means that for every other boundary $(\BD,\lambda)$, up to sets with
measure $0$, there is a measure-preserving surjection of $(\Pi,\nu)$ onto
$(\BD,\lambda)$. A more heuristic interpretation of maximality 
is that $(\Pi,\nu)$ is the finest model for distinguishing limit points at infinity
of the random walk. Existence and uniqueness of the maximal boundary is a general
fact \cite{KaVe}, and it is called the \emph{Poisson boundary}. 

{\sc Kaimanovich~\cite{Kai2}, \cite{KaWo}} has provided a useful geometric tool
for checking maximality. Consider the ``reflected'' right random walk on $\Gamma$
$$
\check Z_n = g_0\gb_1^{-1} \cdots \gb_n^{-1}
$$
starting at $g_0$ (we shall only consider $g_0=\id$).
Its law is the probability measure $\check \mu$ on $\Gamma$, where 
$\check \mu(g) = \mu(g^{-1})$. 

\begin{pro}[{[Kaimanovich].}]\label{pro:strip} 
Let $\mu$ be a probability measure on $\Gamma$ with a finite first moment, 
and let $(\BD,\lambda)$ and $(\check\BD,\check\lambda)$ be a $\mu$- and a 
$\check\mu$-boundary,  respectively. Suppose that there is a measurable 
$G$-equivariant map $S$ assigning to $(\lambda \times \check \lambda)$-almost 
every pair of points  $(\beta,\check \beta)\in \BD\times \check \BD$ 
a non-empty ``strip'' $S(\beta,\check \beta)\subset \Gamma$ such that
for the ball $B(\id,n)$ of radius $n$ in the metric of $\Gamma$,
$$
\frac1n \log \bigl| S(\beta,\check \beta) \cap B(\id,n) \bigr| \to 0 \quad \text{as} \quad
n\to\infty
$$
for $(\lambda \times\check \lambda)$-almost every $(\beta,\check \beta)\in \Pi\times \check \Pi$,
then $(\BD,\lambda)$ and $(\check\BD,\check\lambda)$ are the Poisson boundaries of 
the random walks with law $\mu$ and $\check \mu$, respectively.
\end{pro}

\begin{thm}\label{thm:Poisson} Let $Z_n=(Y_n,X_n)$ be a random walk with law
$\mu$ on $\Gamma =\Zr \wr \Ga$, such that $\mu$ has finite first moment
and $\supp(\mu)$ generates $\Gamma$. If $\Pi$ is defined as in \eqref{eq:Pi} and 
$\nu$ is the limit distribution on $\Pi$ of $(Z_n)$  starting at $\id$, then 
$(\Pi,\nu)$ is the Poisson boundary of the random walk.
\end{thm}

\begin{proof}
By Theorem \ref{thm:conv}, each of the random walks $(Z_n)$ and $(\check Z_n)$
starting at $\id$ converges almost surely to a $\Pi$-valued random variable.
Let $\nu$ and $\check \nu$ be their respective limit distributions.
Then the spaces $(\Pi,\nu)$ and $(\Pi,\check\nu)$ are boundaries of the respective
random walks. If $\beta=(\zeta, \uf) \in \Pi$ then let $\uf(\beta) = \uf$.
Also, if $y \in \T$ and $\uf \in\bd\T$ then then let $\T_y(\uf)$ be the subtree
of $\T$ which is the component of $\uf$ in $\T \setminus \{y\}$ (that is, every ray 
that represents $\uf$ has all but finitely many vertices in this
component of $\T \setminus \{y\}$).

By continuity of $\nu$ and $\check \nu$ (Theorem \ref{thm:conv}), we have
$$
\nu \times \check\nu \bigl(\{ (\beta, \check\beta) \in \Pi \times \Pi : 
\uf(\beta)=\uf(\check\beta) \}\bigr)=0\,.
$$
Therefore, we only need to construct the strips $S(\beta,\check\beta)$ when
$\uf(\beta) \ne \uf(\check\beta)$. 
Thus, let $\beta=(\zeta,\uf), \check\beta = (\check\zeta,\vf) \in \Pi$ and $\uf \ne \vf$.
For any vertex $y$ on the (two-sided infinite) geodesic $\geo{\uf\,\vf}$,
let $\eta_y(\beta,\check\beta)$ be the configuration which coincides with 
$\zeta$ on $\T_y(\vf)$ and with $\check\zeta$ on $\T \setminus \T_y(\vf)$. 
This configuration has finite support, since $\supp(\zeta)$ can only accumulate at
$\uf$ and $\supp(\check\eta)$ can only accumulate at $\vf$, while $\T_y(\vf)$
does not accumulate at $\uf$ and $\T \setminus \T_y(\vf)$ does not accumulate 
at $\vf$. Then define
$$
S(\beta,\check\beta) = \bigl\{ \bigl(\eta_y(\beta,\check\beta),y\bigr) : 
y \in \geo{\uf\,\vf}\bigr\}\,.
$$
This is a subset of $\Gamma$. We check that the map 
$(\beta,\check\beta) \mapsto S(\beta,\check\beta)$
is $\Gamma$-equivariant: let $g = (\eta,x) \in \Gamma$. We have to show that
\begin{equation}\label{eq:equivariant}
g\,S(\beta,\check\beta) = S(g\beta,g\check\beta)\,.
\end{equation}
Now, if $y \in \geo{\uf\,\vf}$ then $xy \in \geo{(x\uf)\,(x\vf)}$. Also,
$T_x\eta_y(\beta,\check\beta) = \eta_{xy}(\beta',\check\beta')$, 
where $\beta' = (T_x\zeta, x\uf)$ and $\check\beta' = (T_x\check\zeta, x\vf)$.
Therefore
$$
\begin{gathered}
\eta + \eta_{xy}(\beta',\check\beta') = \eta_{xy}(\beta'',\check\beta'')\,,
\quad \text{where}\\
\beta'' = (\eta + T_x\zeta, x\uf) = g\beta 
\AND \check\beta'' = (\eta + T_x\check\zeta, x\vf) = g\check\beta\,.
\end{gathered}
$$
Thus, for $y \in \geo{\uf\,\vf}$,
$$
g\bigl(\eta_y(\beta,\check\beta),y\bigr) = 
\bigl( \eta+T_x\eta_y(\beta,\check\beta) , xy \bigr) 
= \bigl( \eta_{xy}(g\beta,g\check\beta) , xy \bigr) 
$$
This proves \eqref{eq:equivariant}.

Finally, if $\eta\,y \in S(\beta,\check\beta)$ and $d(\zero\,o,\eta\,y) \le n$
then $d(o,y) \le n$. Since 
$$
|\{ y \in \geo{\uf\,\vf} : d(o,y) \le n \}| \le 2n\,,
$$
we see that all conditions of Proposition \ref{pro:strip} are satisfied.
\end{proof}

\section{The Dirichlet problem at infinity}\label{Dirichlet}

In this section we shall assume in addtion that our random walk
on $\Gamma$ is \emph{irreducible} in the sense that its law $\mu$
is such that $\supp(\mu)$ generates $\Gamma$ as \emph{semi}group. 
Equivalently, this means that for every pair of elements $g= \eta\,x$,
$g' = \eta'\,x' \in \Gamma$ the probability that the random walk
starting at $g$ ever visits $g'$ is $>0$. 

Also, it will be convenient to consider the limit distribution $\nu$ 
of \eqref{eq:nu} as a Borel measure on the compact set 
$\overline{\Pi} = \whCC \times \bd\T\,.$
The irreducibility hypothesis implies that
\begin{equation}\label{eq:suppnu}
\supp(\nu) = \overline{\Pi}\,,
\end{equation}
that is, the whole of $\overline{\Pi}$ is active (as we shall see below).

With respect to our random walk with transition 
probabilities $p(g,g')=\mu(g^{-1}g')$ on $\Gamma$, a function
$h:\Gamma \to \R$ is called \emph{harmonic,} if it satisfies the 
weighted mean value property
$$
h(g) = \sum_{g'}p(g,g')h(g') \quad\text{for all}\; g \in \Gamma\,.
$$ 
For $g \in \Gamma$, we define $\nu_g$ by $\nu_g(U)=\nu(g^{-1} U)$, where 
$U$ runs through Borel subsets of $\overline{\Pi}$.  
It is a basic feature of the Poisson boundary that every bounded harmonic
function $h$ on $\Gamma$ has a unique integral representation
\begin{equation}\label{eq:poissint}
h(g) = \int_{\whCC \times \bd\T} f \,d\nu_g\,,
\end{equation}
where $f \in L^{\infty}(\nu)$, see e.g. \cite{KaVe}.
 
Conversely, once we have convergence of the random walk to the boundary, 
any function $f \in L^{\infty}(\nu)$ gives rise to a harmonic
function via the integral formula \eqref{eq:poissint}. In fact, for this we do
not need to know that $\overline{\Pi}$ is the Poisson boundary; the only point 
is that otherwise we will not get \emph{all} bounded harmonic functions
via  \eqref{eq:poissint}.

Related to existence of the limit measure on the boundary, there is the 
question whether the \emph{Dirichlet problem at infinity} is solvable. 
In our case it reads as follows:

Does every continuous function on $\overline{\Pi}=\whCC \times \bd\T$ extend 
continuously to a function on $\Gamma \cup \overline{\Pi}$
which is harmonic on $\Gamma$\,?

If the answer is positive, then that harmonic function must be given by 
\eqref{eq:poissint}, and 
we would like that whenever $f$ is continuous it should hold that
$\lim_{g\rightarrow \beta}h(g)=f(\beta)$
for every $\beta\in \overline{\Pi}$. 
We then say that the \emph{Dirichlet problem at infinity is solvable}
for continuous functions on $\overline{\Pi}$.

We note that usually, the Dirichlet problem at infinity refers to the
harmonic extension of  any continuous function that is given on the 
\emph{whole} boundary in a (suitable) compactification of the state space, 
compare with \cite[\S 20]{Wbook}. In our case, the whole boundary is the set 
$(\wh\CC \times \wh\T) \setminus (\CC \times \T)$
which contains $\whCC \times \bd\T$ as a proper, compact subset. 
However, the complement of $\whCC \times \bd\T$ is not charged by $\nu$,
so that boundary data given on that complement have no effect on the
harmonic function $h$ of \eqref{eq:poissint}, and we cannot expect continuity
at those points. Therefore we have to restrict to continuous functions on
$\whCC \times \bd\T$.

The \emph{Green kernel} of the projected random walk $(X_n)$ on 
$\T \equiv \Ga$ with law $\wt \mu$ as in \eqref{eq:mutilde} is
$$
\wt G(x,y) = \sum_{n=0}^{\infty} \wt p^{(n)}(x,y) = 
\sum_{n=0}^{\infty} \wt\mu^{(n)}(x^{-1}y)\,,\quad x, y \in \T\,,
$$
where $\wt p^{(n)}(\cdot,\cdot)$ denotes $n$-step transition probabilities
and $\wt\mu^{(n)}$ is the $n$-th convolution power of $\wt \mu$. By
irreducibility and transience, $0 < \wt G(x,y) < \infty$. This is the 
expected number of times that $(X_n)$ visits $y$, given that $X_0=0$. 

\begin{lem}\label{lem:green}
If $\wt \mu$ is irreducible, then Green kernel vanishes at infinity, that is,
$$
\lim_{d(x,y) \to \infty} \wt G(x,y) = 0\,.
$$
\end{lem}

\begin{proof} The Group $\Ga$ is non-amenable. Let $S = \supp(\wt \mu)$.
Since $S$ generates $\Ga$ as a semigroup, it is a well-known exercise to show
that the sugroup of $\Ga$ generated by $S\cdot S^{-1}$ is a finite-index
normal subgroup, whence also non-amenable. Therefore one can apply 
Th\'eor\`eme 2 of {\sc Derriennic and Guivarc'h} \cite{DG} and/or
Th\'eor\`eme 2 of {\sc Berg and Christensen} \cite{DG} to obtain that
the measure $\sum_{n=0}^{\infty} \wt\mu^{(n)}$ defines a bounded convolution
operator on $\ell^2(\Ga)$. It follows that the Green kernel vanishes
at infinity. This may also be deduced by applying the main theorem of
\cite{Ke2} to $\wt \mu * \check{\wt \mu}$. 
\end{proof}

Below, we shall need the quantity 
$$
\wt F(x,y) = \wt G(x,y)/\wt G(y,y)
= \Prob[\exists\ n \ge 0: X_n = y \mid X_0 = 0]\,,
$$
which also vanishes at infinity.

\begin{thm}\label{thm:Diri}
If $\mu$ is irreducible and the lamplighter operates at finite range,
then the Dirichlet problem at infinity for continuous functions on 
$\overline{\Pi} = \whCC \times \bd\T$ is solvable.
\end{thm}

\begin{proof}
The typical ``probabilistic'' method for proving this (see e.g. \cite[\S 20]{Wbook}) 
is to show the following: (i) the random walk $Z_{n}$ converges to the boundary, 
and (ii) for the corresponding harmonic measure class 
$\{\nu_g : g=\eta\, x \in \Gamma\}$, one has
\begin{equation}\label{eq:weakly}
\lim_{g\rightarrow \beta}\nu_{g}=\delta_{\beta}
\quad \text{weakly for every}\; \beta = (\zeta,\uf) \in \whCC \times \bd\T\,.
\end{equation}
Point (i) is affirmed by Theorem \ref{thm:conv}. For proving (ii),
it will be convenient to consider $\nu_g$ as a measure on 
$\whCC \times \wh\T$ which charges only the set 
$\overline{\Pi}= \whCC \times \bd\T$.
Thus, we show that for any neigbourhood $U$ in $\whCC \times \wh\T$ of  
$\beta = (\zeta,\uf) \in \overline{\Pi}$, we have 
$$
\lim_{g \to \beta} \nu_g(U^c)=0\,,
$$
where $U^c= (\whCC \times \wh\T)\setminus U$. Here, it is sufficient to 
take $U$ in a suitable neighbourhood basis of $\beta$. Such a basis is obtained
as follows: take a vertex $y \in \T$ and let 
$\wh\T_y(\uf)$ be the closure in $\wh\T$ of the subtree $\T_y(\uf)$ (the
component of $\uf$ in $\T \setminus \{y\}$, see the proof of Theorem
\ref{thm:Poisson}). The familiy of all $\wh T_y(\uf)$, 
$y \in \T$, is a neighbourhood
basis of $\uf$ in $\wh \T$. Also, the family of all sets
$\wh C_A(\zeta) = \{ \zeta' \in \whCC : \zeta' = \zeta \;\text{on}\;A \}$,
where $A \subset \T$ is finite, is a neigbouhrhood basis of $\zeta$
in $\whCC$ for the topology of pointwise convergence of configurations.
Thus, 
$$
\{ U_{y,A}(\beta)=
\wh T_y(\uf) \times \wh C_A(\zeta) : y \in \geo{o\,\uf}\,,\; A \subset \T 
\;\text{finite} \}
$$
is a neighbourhood basis of $\beta=(\zeta,\uf)$. Consider $U = U_{x,A}(\beta)$.
Then 
$$
U^c \subset V_1 \cup V_2\,,\quad\text{where}\quad
V_1 = \bigl(\wh\T \setminus \wh\T_y(\uf)\bigr) \times \whCC \AND V_2 = 
\wh \T_y(\uf) \times \bigl(\whCC \setminus \wh C_A(\zeta)\bigr)\,.
$$
We now prove that $\nu_g(V_i) \to 0$ for $i=1,2$ as $g=\eta\,x \to \beta$.

Regarding $V_1$, we first remark that the Dirichlet problem at infinity
for continuous functions on $\bd\T$ is solvable for the random walk $(X_n)$
with law $\wt\mu$ on $\Ga \equiv \T$, see \cite{CaSo} and 
\cite[Cor. 21.12]{Wbook}. Therefore
$$
\lim_{g \to \beta} \nu_g(V_1) 
= \lim_{x \to \uf} \Prob_x[X_{\infty} \in \bd\T \setminus \wh T_y(\uf)] = 0\,.
$$
Regarding $V_2$, let $R = R(\mu)$ be the bound on the range of 
\eqref{eq:lamprange}. Set $A^R = \{v \in \T : d(v,A) \le R \}$.
Suppose $\eta = \eta(g)$ is sufficiently close to the limit
$\zeta$ so that $\eta(v)=\zeta(v)$ for every $v \in A$. If the random walk
$Z_n=(Y_n,X_n)$ starting at $g$ converges to a limit point in $V_2$
then $(X_n)$ must visit $A^R$ in order to modify the states of the 
lamps at the points in $A$. Therefore
$$
\nu_g(V_2) \le \sum_{v \in A^R} \wt F(x,v)\,,
$$ 
which tends to zero when $x \to \uf$ by Lemma \ref{lem:green}.
This concludes the proof.
\end{proof}

\section{Final remarks}\label{sec:remarks}

Theorem \ref{thm:Poisson} is yet one more application of the very useful
strip criterion of {\sc V. A. Kaimanovich} who, in 
private comunication, 
has informed us that in the unpublished paper \cite{Kai3} he uses a method in a 
somewhat similar spirit to the proof of Theorem \ref{thm:Poisson} to describe 
the Poisson boundary for 
random walks on $\Zr \wr \Z^d$, where the projection of the random walk onto 
the integer grid $\Z^d$ has non-zero drift.

\medskip

Regarding the Dirichlet problem at infinity, we repeat here that this type 
of question can be asked whenever one has a compactification of the state
space in whose topology the random walk converges almost surely to
the boundary at infinity. Solvability of the Dirichlet problem for
continuous functions on the boundary (or rather its active part, i.e., 
the support
of the limit measure) is by no means the same as having determined
the Poisson (or even Martin) boundary, as one finds erroneously stated every 
now and then in published work. For example, in our case, we know that 
$(\overline{\Pi},\nu)$ is the Poisson boundary when $\mu$ has finite first
moment, while we proved that the corresponding Dirichlet problem is
solvable when the lamplighter operates at finite range, but $\mu$ need not
have finite first moment for that. Thus, we have given positive answers
to both questions simultaneously only when $\mu$ has finite first moment
\emph{and} the lamplighter operates at finite range. 

\medskip

The methods that we have used can be extended to space-homogenous (in the sense
of \cite{KaWo}) lamplighter random walks over hyperbolic graphs, graphs with 
infinitely many ends, and other classes of transitive base graphs that can  
be used in  the construction of lamplighter graphs according to \eqref{eq:nbhd}.
Also, the fact that the states of the lamps are encoded by the complete graph
$\Kr$, resp. the cyclic group $\Zr$, does not play an essential role. 
Finally, with additional effort, the result regarding the Dirichlet problem 
can apparently be extended.
The detailed elaboration of these general facts remains reserved to future
work, while the present note has aimed at giving a short and hopefully readable 
explanation of the basic aspects.


\begin{thebibliography}{22}

\bibitem{BaWo} Bartholdi, L., and Woess, W.: \emph{Spectral computations on 
lamplighter groups and Diestel-Leader graphs,}  J. Fourier Anal. Appl.
{\bf 11} (2005) 175 - 202.

\bibitem{BC} Berg, C., and Christensen, J. P. R.: \emph{Sur la norme des
op\'erateurs de convolution,} Invent. Math. {\bf 23} (1974) 173--178.

\bibitem{BrWo1} Brofferio, S., and Woess, W.: 
\emph{Green kernel estimates and the full Martin boundary for random walks
on lamplighter groups and Diestel-Leader graphs}, Annales Inst. H. Poincar\'e 
(Prob. \& Stat.)  {\bf 41} (2005) 1101--1123.  

\bibitem{BrWo2} Brofferio, S., and Woess, W.: 
\emph{Positive harmonic functions for semi-isotropic random walks on trees, 
lamplighter groups, and DL-graphs}, Potential Analysis {\bf 24} (2006) 245--265.

\bibitem{CaKaWo} Cartwright, D. I., Kaimanovich, V. A., and Woess, W.: 
\emph{Random walks on the affine group of local fields and of homogeneous 
trees,} Ann. Inst Fourier (Grenoble) {\bf 44} (1994) 1243--1288.

\bibitem{CaSo} Cartwright, D. I., and Soardi, P. M.: \emph{Convergence
to ends for random walks on the automorphism group of a tree,} 
Proc. Amer. Math. Soc. {\bf 107} (1989) 817--823 

\bibitem{CeWo} Ceccherini-Silberstein, T., and Woess, W.: \emph{Growth and
ergodicity of context-free languages}, 
Trans. Amer. Math. Soc. {\bf 354} (2002) 4597--4625.

\bibitem{DG} Derriennic, Y., and Guivarc'h, Y.: \emph{Th\'eor\`eme de 
renouvellement pour les groupes non moyennables,} C. R. Acad. Sci. Paris,
S\'erie A {\bf 277} (1973) 613--615.
  
\bibitem{De} Derriennic, Y.: \emph{Quelques applications du th\'eor\`eme ergodique 
sous-additif,} Ast\'erisque {\bf 74} (1980) 183--201.

\bibitem{DiSc} Dicks, W., and Schick, Th.: \emph{The spectral measure 
of certain elements of the complex group ring of a wreath product,}  
Geom. Dedicata  {\bf 93}  (2002) 121--137.

\bibitem{Er1} Erschler, A. G.: \emph{On the asymptotics of the rate 
of departure to infinity} (Russian),  Zap. Nauchn. Sem. S.-Peterburg. 
Otdel. Mat. Inst. Steklov. (POMI) {\bf 283} (2001) 251--257, 263. 

\bibitem{Er2} Erschler, A. G.: \emph{On drift and entropy growth for
random walks on groups,}  Ann. Probab. {\bf 31} (2003) 1193--1204.

\bibitem{Er3} Erschler, A. G.: \emph{Liouville property for groups and
manifolds,} Invent. Math. {\bf 155} (2004) 55--80. 

\bibitem{Fu} Furstenberg, H.: \emph{Random walks and discrete 
subgroups of Lie groups.} In \emph{Advances in Probability and Related Topics,}
{\bf 1} (P. Ney, ed.), pp. 1--63, M. Dekker, New York,  1971.

\bibitem{Gi} Gilch, L. A.: \emph{Rate of escape on the lamplighter tree,}  
Preprint, TU Graz (2006).

\bibitem{GrZu} Grigorchuk, R. I., and \.Zuk, A.: \emph{The lamplighter 
group as a group generated by a 2-state automaton, and its spectrum,}  
Geom. Dedicata {\bf 87} (2001) 209--244. 

\bibitem{Kai1} Kaimanovich, V. A.: \emph{Poisson boundaries of random 
walks on discrete solvable groups,} in: \emph{Probability Measures on Groups X}
(ed. H. Heyer), pp. 205--238,  Plenum, New York, 1991.

\bibitem{Kai2} Kaimanovich, V. A.: \emph{The Poisson formula for groups 
with hyperbolic properties,} Annals of Math. {\bf 152} (2000) 659--692.

\bibitem{Kai3} Kaimanovich, V. A.: \emph{Poisson boundary of discrete groups, a survey,} 
unpublished manuscript, {\tt http://name.math.univ-rennes1.fr/vadim.kaimanovich/list.htm}

\bibitem{KaVe} Kaimanovich, V. A., and Vershik, A. M.: 
\emph{Random walks on discrete groups: boundary and entropy,} Ann.
Probab. {\bf 11} (1983) 457--490.

\bibitem{KaWo} Kaimanovich, V. A., and Woess, W.: \emph{Boundary and
entropy of space homogeneous Markov Chains,} Ann. Probab. {\bf 30} (2002) 
323-363.

\bibitem{Ke1} Kesten, H.: \emph{Symmetric random walks on groups,}
Trans. Amer. Math. Soc. {\bf 92} (1959) 336--354.

\bibitem{Ke2} Kesten, H.: \emph{Full Banach mean values on countable groups,}
Math. Scand. {\bf 7} (1959) 146--156.

\bibitem{LyoPemPer} Lyons, R., Pemantle, R., and Peres, Y.: 
\emph{Random walks on the lamplighter group,} Ann. Probab. {\bf 24} (1996)
1993--2006.

\bibitem{Pa} Parry, W.: \emph{Growth series of some wreath products}, 
Trans. Amer. Math. Soc. {\bf 331} (1992) 751--759.

\bibitem{PitSal1} Pittet, C., and Saloff-Coste, L.: \emph{Amenable groups, 
isoperimetric profiles and random walks,} in: \emph{Geometric Group 
Theory Down Under} (Canberra, 1996),  pp. 293--316, de Gruyter, Berlin, 1999. 

\bibitem{PitSal2} Pittet, C., and Saloff-Coste, L.: \emph{On random walks 
on wreath products,}  Ann. Probab. {\bf 30} (2002) 948--977.

\bibitem{Rev1} Revelle, D.: \emph{Rate of escape of random walks 
on wreath products,} Ann. Probab. {\bf 31} (2003) 1917--1934.

\bibitem{Rev2} Revelle, D.: \emph{Heat kernel asymptotics on the 
lamplighter group,} Electron. Comm. Probab. {\bf 8} (2003), 142--154.

\bibitem{Wbook} Woess,  W.: \emph{Random Walks on Infinite Graphs and Groups}, 
Cambridge Tracts in Mathematics {\bf 138}, Cambridge University Press,
Cambridge, 2000.

\bibitem{Wo1} Woess,  W.: \emph{Lamplighters, Diestel-Leader graphs, 
random walks, and harmonic functions,} Combinatorics, Probability \& Computing
{\bf 14} (2005) 415--433.

\bibitem{Wo2} Woess,  W.: \emph{A note on the norms of transition operators 
on lamplighter graphs and groups,}  Internat. J. Algebra Comput.
{\bf 15} (2005) 1261--1272.

\end{thebibliography}
\end{document}